\documentclass[10pt,letterpaper,english]{article}
\usepackage[T1]{fontenc}
\usepackage[latin1]{inputenc}
\usepackage{amsmath}
\usepackage{graphicx}
\usepackage{amssymb}
\makeatletter

\makeatletter
\usepackage{amscd}

\usepackage{amsfonts}

\usepackage{babel}\makeatother

\begin{document}
\newtheorem{proposition}{Proposition}[section] \newtheorem{definition}{Definition}[section]
\newtheorem{corollary}{Corollary}[section] \newtheorem{lemma}{Lemma}[section]
\newtheorem{theorem}{Theorem}[section] \newtheorem{example}{Example}[section]

\title{\textbf{An Optimal Control Perspective on \\ Classical
		and Quantum Physical Systems.}}

\author{M. Contreras G.\thanks{Universidad Metropolitana de Ciencias de la Educación UMCE, Chile, email: mauricio.contreras@umce.cl} \ and \ M. Villena\thanks{Universidad Adolfo Ib\'añez, Chile, email: marcelo.villena@uai.cl} }
\maketitle
\noindent

\begin{abstract}
	In this paper, we analyze classical and quantum physical systems from an optimal control perspective. Specifically, we explore whether their associated dynamics can correspond to an open or closed-loop feedback evolution of a control problem.
	Firstly, for the classical regime, when it is viewed in terms of the theory of canonical transformations, we find that it can be described by a closed-loop feedback problem. Secondly, for a quantum physical system, if one realizes that the Heisenberg commutation relations themselves can be thought of as constraints in a non-commutative space, then the momentum must be dependent on the position for any generic wave function. This implies the existence of a closed-loop strategy for the quantum case. Thus, closed-loop feedback is a natural phenomenon in the physical world. For the sake of exposition, we give a short review of control theory, and some familiar examples at the classical and quantum levels are analyzed.
\end{abstract}

\maketitle

\section{Introduction}

In the last two decades, there has been increasing interest among physicists in applying ideas from physics to finance and economics, as one can see in some classical texts.\cite{stanley} \ \cite{boucheaud} \ \cite{baaquie} In this paper, we want to do the exact opposite, that is, we apply ideas from the optimal control theory usually used in finance and
economics,\cite{kamien} \  \cite{sheti} \ \cite{caputo} \ \cite{weitzman} \ \cite{dockner} to the classical and quantum physical world. Indeed, open and closed-loop control problems are not in the toolbox of mathematical methods in physics; nevertheless, dynamic optimization with its optimal control theory is the corner-stone of modern economic analysis. \\

There is not much literature associated with this endeavor.
Recently, there have appeared some studies of how a control problem associated with an economic model can be interpreted as a second class constrained physical system.\cite{contrerasRPMV} \ \cite{hojman} \ \cite{teturo} \ \cite{contreras2} In [9] it is found that at the classical level, the constrained dynamics given by Dirac's brackets are the same as the dynamics given by the Pontryagin equations \cite{Pontryagin}. The right quantization of this second-class constrained system ended with a Schrödinger equation
that is just the Hamilton--Jacobi--Bellman equation used in optimal control theory and economics \cite{Bellman}. It was surprising to find that after the quantization of the Pontryagin theorem, Bellman's maximum principle is obtained. \\

From the perspective of control theory,  this is an exciting finding in terms of the feedback of these systems. The Pontryagin
theory is characterized by open-loop strategies for the Lagrangian multiplier, whereas the Hamilton--Jacobi--Bellman equation has closed-loop strategies intrinsically. Since Pontryagin's theory is equivalent to a classical mechanical model and Bellman's theory is just a quantum mechanical one,
this induces one to relate open-loop strategies with classical mechanics and closed-loop strategies with quantum mechanics. All this sounds as if the quantizations of open-loop strategies would be the closed-loop ones. \\

Thus, if an economic feedback system characterized by an optimal control theory model can be seen as a physical system, one can naturally ask the inverse question: can the physical systems be seen from an optimal control perspective? Or, more specifically, do there exist open-loop and closed-loop strategies in physical systems?  If the answer is affirmative, then, where are they? Moreover, how and why do these feedback problems appear? This paper gives some clues and answers to all these questions. \\

For this purpose, we found that the classical theory must be represented as a closed-loop feed-back problem. Indeed, due to the invariance of the classical Hamiltonian equations of motion (which represent an open-loop dynamics) by canonical transformations, both schemes are equivalent. In the same way, from the perspective of quantum physics, we see that if one identifies the momentum with the Lagrange multiplier, the relation between momentum and position is just a closed-loop strategy in optimal control theory. Thus, both open and closed-loop feed-backs are also a natural phenomenon in our physical world. \\

In the following three sections we review, for the non specialist reader, the ideas of closed and open-loop strategies that appear in control theory.

\section{Dynamic Optimization: The Pontryagin Approach}

To emphasize the fundamental ideas and to keep the equations simple, a generic one-dimensional physical system will be considered. Generalizations to higher dimensions are straightforward. Consider an optimal control problem commonly used in financial applications.\cite{kamien} \ \cite{sheti} It is required
to optimize the cost functional
\begin{equation}
A[x,u]=\int_{t_{0}}^{t_{1}}F(x,u,t)\ dt,
\end{equation}
where $x$ represents a state variable (for example, the production of a certain article) and $u$ is a control variable (such as the marketing costs). The state variable must satisfy the market dynamics
\begin{equation}\label{xdotf}
\dot{x}=f(x,u,t)   \ \ \ \ \ \  	x(t_{0})=x_{0}.
\end{equation}
The problem is to determine how to obtain the production trajectory $x=x(t)$ and the control path $u=u(t)$ to optimize the cost functional.
To get the solution, the method of Lagrange multipliers is applied, so an improved functional $A$ on the extended configuration space $(x,u,\lambda)$
is considered, which is defined by
\begin{equation}
A[x,u,\lambda]=\int_{t_{0}}^{t_{1}}F(x,u,t)-\lambda(\dot{x}-f(x,u,t))\ dt.\label{A}
\end{equation}
To obtain the solution, the integrand of (\ref{A}) can be interpreted
as the Lagrangian:
\begin{equation}
L(x,u,\lambda,\dot{x},\dot{u},\dot{\lambda})=F(x,u,t)-\lambda(\dot{x}-f(x,u,t)).\label{Lagrangiancontopt1}
\end{equation}
The extremal curves then satisfy the Euler--Lagrange equations:
\begin{equation}
\frac{\partial L}{\partial\lambda}-\frac{d}{dt}\ \big(\frac{\partial L}{\partial\dot{\lambda}}\big)=0 , \ \ \ \ \ \ \ \frac{\partial L}{\partial x}-\frac{d}{dt}\ \big(\frac{\partial L}{\partial\dot{x}}\big)=0, \ \ \ \ \ \ \  \frac{\partial L}{\partial u}-\frac{d}{dt}\ \big(\frac{\partial L}{\partial\dot{u}}\big)=0.
\end{equation}
These are also written respectively as
\begin{equation} \label{pontry}
\dot{x}=\frac{\partial H}{\partial\lambda} ,   \ \ \ \ \ \ \ \ \ \ \ \dot{\lambda}=-\frac{\partial H}{\partial x} ,  \ \ \ \ \ \ \ \ \ \ \
\frac{\partial H}{\partial u}=0 ,
\end{equation}
with $H$ defined by
\begin{equation}
H=H(x,u,\lambda)=F(x,u,t)+\lambda f(x,u,t).
\end{equation}
Equations (\ref{pontry}) are the well-known Pontryagin equations, which are obtained in optimal control theory, through the Pontryagin maximum principle, and give the solution to the optimization problem.

\section{Open-loop and closed-loop strategies}

The action (\ref{A}) can be written in a compact form as
\begin{equation}
A[x,u,\lambda]=\int_{t_{0}}^{t_{1}}-\lambda\dot{x}+H(x,u,\lambda,t)\ dt.\label{A2}
\end{equation}
Strictly speaking, the Pontryagin equations must be obtained by optimizing the action (\ref{A2}) with respect to its three variables $x,u,\lambda$. If $\delta x,\ \delta u$ and $\delta\lambda$ are the corresponding functional variations of the initial variables, by expanding the
Hamiltonian in a Taylor series, keeping the first-order terms only, integrating by parts, and using the fact that $\delta x(t_{0})=0$ (the initial point is fixed by the initial condition (\ref{xdotf})), then one obtains
\begin{equation}
\delta A={\displaystyle \int_{t_{0}}^{t_{1}}\Bigl[\left(\frac{\partial H}{\partial\lambda}-\dot{x}\right)\delta\lambda+\left(\frac{\partial H}{\partial x}+\dot{\lambda}\right)\delta x+\frac{\partial H}{\partial u}\delta u\Bigr]dt-\lambda(t_{1})\delta x(t_{1}).}\label{A3}
\end{equation}
To maximize the action, all the first-order terms in $\delta x,\delta u$
and $\delta\lambda$ must vanish.  Now, it is well-known that there are two classes of control strategies: 
\begin{itemize}
	\item  \textbf{open-loop strategies} that depend only on time: $u=u(t)$,
	and 
	\item  \textbf{closed-loop strategies} that depend on the state variable
	$x$ and time: $u=u(x,t)$. \cite{erikson1} 
\end{itemize} 

In the case of an open-loop strategy, the variables $x=x(t), u=u(t)$ and $\lambda = \lambda(t)$ are independent, and so $\delta x,\delta u$ and $\delta\lambda$ are linearly independent. Hence, the Pontryagin equations (\ref{pontry}) and the transversality condition $\lambda(t_{1})=0  $ can be obtained from Equation (\ref{A3}). In fact, the solutions of the
Pontryagin equations give the extremal curves
\begin{equation} \label{optimaluaster}
x=x^{*}(t) , \ \ \ \ \ \ \ \ \ \ \ \  \lambda=\lambda^{*}(t) ,
\ \ \ \ \ \ \ \ \ \ \ \ 
u=u^{*}(t) ,
\end{equation} 
which are the solutions of Equations (\ref{pontry}). But what happens, however, with closed-loop strategies? Here,
due to the relations between the variables in $u=u(x,t)$, the functional variations are related by $\delta u=\frac{\partial u}{\partial x}\delta x$. Substituting this into (\ref{A3}) yields
\begin{equation}
\delta A={\displaystyle \int_{t_{0}}^{t_{1}}\Bigl[\left(\frac{\partial H}{\partial\lambda}-\dot{x}\right)\delta\lambda+\left(\frac{\partial H}{\partial x}+\frac{\partial H}{\partial u}\frac{\partial u}{\partial x}+\dot{\lambda}\right)\delta x\Bigr]dt-\lambda(t_{1})\delta x(t_{1}).}
\end{equation}
If $\lambda$ and $x$ remain independent, the optimization of the functional implies
\begin{equation}
\frac{\partial H}{\partial\lambda}-\dot{x}=0\;,\qquad\quad\frac{\partial H}{\partial x}+\frac{\partial H}{\partial u}\frac{\partial u}{\partial x}+\dot{\lambda}=0,
\end{equation}
plus the transversality condition, but the equation that gives the optimal condition for the control in (\ref{pontry}) is lost. Then, if $u$ is not known as a fixed function of $x$ from the beginning, that is, one knows that the constraints $u = u(x, t)$ exist, but one does not know their specific form, one has three unknowns $x,u$, and $\lambda$ but only two equations of motion. Thus, the variational problem is inconsistent for an arbitrary closed-loop strategy $u=u(x,t)$ because the functional form of $u$ is not determined by the equations of motion. It must be given from the beginning. \\
Now, note that for a open-loop strategy, the control equation in (\ref{pontry}) is, in fact, the following algebraic equation for $u$:
\begin{equation}
\frac{\partial H}{\partial u} =	\frac{\partial F(x,u,t)}{\partial u}+\lambda\ \frac{\partial f(x,u,t)}{\partial u}=0.\label{pontryuv2}
\end{equation}
From this equation, the optimal control $u^*$ can be obtained in principle, as a function of $x$ and $\lambda$:
\begin{equation}
u=u^{*}(x,\lambda,t),\label{uxl}
\end{equation}
One can conclude, then, that the same optimization problem implies that the naive optimal open-loop $u^{*}$ strategy in (\ref{optimaluaster}) is actually just a closed-loop strategy! How can this be consistent with the variational problem (\ref{A3}), in which $x$, $u$ and $\lambda$
are independent variables? \\
To understand this, consider an arbitrary closed-loop strategy of the general form $u=u(x,\lambda,t)$.  Then \ $\delta u=\frac{\partial u}{\partial x}\delta x+\frac{\partial u}{\partial\lambda}\delta\lambda$.
After substituting this into (\ref{A3}), $\delta A$ equals
\begin{equation}
{\displaystyle \int_{t_{0}}^{t_{1}}\Bigl[\left({\displaystyle \frac{\partial H}{\partial\lambda}+{\displaystyle \frac{\partial H}{\partial u}\frac{\partial u}{\partial\lambda}-\dot{x}}}\right)\delta\lambda+\left(\frac{\partial H}{\partial x}+\frac{\partial H}{\partial u}\frac{\partial u}{\partial x}+\dot{\lambda}\right)\delta x\Bigr]dt-\lambda(t_{1})\delta x(t_{1}).}\label{A3v2}
\end{equation}
If $x(t)$ and $\lambda(t)$ are independent variables, we have the following equations of motion from the variational principle (\ref{A3v2}):
\begin{equation} \label{ponl2}
\frac{\partial H}{\partial\lambda}+\frac{\partial H}{\partial u}\frac{\partial u}{\partial\lambda}-\dot{x}=0,  \ \ \ \ \ \ \ \ \ \  \frac{\partial H}{\partial x}+\frac{\partial H}{\partial u}\frac{\partial u}{\partial x}+\dot{\lambda}=0.
\end{equation}
Obviously these equations, for an arbitrary closed-loop strategy
$u=u(x,\lambda,t)$, differ from the Pontryagin open-loop equations. Now, choose $u$ as the special optimal closed-loop strategy $u^{*}$ which is the solution of (\ref{pontryuv2}), so Equations (\ref{ponl2}) are then reduced to
\begin{equation} \label{ll1}
\frac{\partial H(x,u^{*}(x,\lambda,t),\lambda,t)}{\partial\lambda}-\dot{x}=0, \ \ \ \ \ \ \ \ \ \ \frac{\partial H(x,u^{*}(x,\lambda,t),\lambda,t)}{\partial x}+\dot{\lambda}=0.
\end{equation}
But this last set of two equations is equivalent to the three Pontryagin equations. The solutions of (\ref{ll1}) give optimal
paths $x=x^{*}(t)$ and $\lambda=\lambda^{*}(t)$, from which the
optimal control open-loop strategy $u^{*}(t)$ in (\ref{optimaluaster}) can be computed given the optimal closed-loop strategy (\ref{uxl}), by
\begin{equation} \label{openloop_u*(t)}
u=u^{*}(t)=u^{*}(x^{*}(t),\lambda^{*}(t),t)
\end{equation}
Thus, the special optimal closed-loop strategy $u=u^{*}(x,\lambda,t)$ is completely equivalent to an open-loop strategy $u=u^{*}(t)$ over the optimal $x^{*}(t)$ and $\lambda^{*}(t)$ trajectories. Then, the optimal closed-loop strategy $u^*(x(t),\lambda(t), t)$ is the same object given by the open-loop one from a dynamical point of view. So from now on, we will not distinguish between them for the case of Pontryagin's theory.\\
Note that for an arbitrary given closed-loop strategy $u = u(x, \lambda,t)$ which is not optimal ($\frac{\partial H}{\partial u} \ne 0$), Equations (\ref{ponl2}) mean that the action is optimized, but these extremals are not necessarily global. In fact, its dynamics is not given by Pontryagin's equations. The optimal condition $\frac{\partial H}{\partial u} = 0$ in (\ref{pontry}) gives the global extremal for the action $A$, as was established by Pontryagin's maximum principle.\\

To end this section, suppose that $\lambda$ and $x$ are not independent but are related by $\lambda=\lambda(x,t)$; then, the variation of $\lambda$ is  $\delta\lambda=\frac{\partial\lambda}{\partial x}\ \delta x    $.
Substituting this into (\ref{A3v2}), choosing the optimal control strategy $u=u^{*}$, using the transversality condition $\lambda(t_{1})=0$
and the fact that $\dot{\lambda}={\displaystyle \frac{\partial\lambda}{\partial x}\dot{x}+\frac{\partial\lambda}{\partial t}}$,
we get
\begin{equation}
\delta A=\int_{t_{0}}^{t_{1}}\Bigl[\left(\frac{\partial H}{\partial\lambda}\frac{\partial\lambda}{\partial x}-\dot{x}\frac{\partial\lambda}{\partial x}\right)+\left(\frac{\partial H}{\partial x}+\frac{\partial\lambda}{\partial x}\dot{x}+\frac{\partial\lambda}{\partial t}\right)\Bigr]\delta x\ \ dt,\label{A3v3}
\end{equation}
but, as $H=F(x,u^{*},t)+\lambda f(x,u^{*},t)$, it follows that 
\[
\frac{\partial H}{\partial\lambda}=f(x,u^{*},t,) \ \ \ \ \ \ \ \ \ \ \ \ \ \ \ \  \frac{\partial H}{\partial x}=\frac{\partial F(x,u^{*},t)}{\partial x}+\lambda\frac{\partial f(x,u^{*},t)}{\partial x}.
\]
Thus, (\ref{A3v3}) gives finally
\[
\delta A=\int_{t_{0}}^{t_{1}}\Bigl[\left(f\frac{\partial\lambda}{\partial x}+\frac{\partial F}{\partial x}+\lambda\frac{\partial f}{\partial x}\right)+\frac{\partial\lambda}{\partial t}\Bigr]\delta x\ \ dt
\]
or
\[
\delta A=\int_{t_{0}}^{t_{1}}\Bigl[\frac{dH^{*}}{dx}+\frac{\partial\lambda}{\partial t}\Bigr]\delta x\ \ dt.
\]
where
\begin{equation}
H^{*}=H^{*}(x,t)=H(x,u^{*}(x,\lambda(x,t),t),\lambda(x,t),t)\label{A3v66}
\end{equation}
is the reduced Hamiltonian in terms of $x$. In this way, the optimization
of the action $A$ implies that the closed-loop $\lambda=\lambda(x,t)$
strategy must satisfy the optimal consistency condition
\begin{equation}
\frac{dH^{*}(x,t)}{dx}+\frac{\partial\lambda}{\partial t}=0.\label{A3v6}
\end{equation}
Equation (\ref{A3v6}) is closely related to the Hamilton--Jacobi--Bellman
equation. In fact, if the Lagrangian multiplier is of the form $\lambda(x,t)=\frac{\partial J(x,t)}{\partial x}$,
then Equation (\ref{A3v6}) implies that $J(x,t)$ satisfies the Hamilton--Jacobi--Bellman equation.\cite{contrerasRPMV} \\
Let $\lambda^{*}(x,t)$ be a solution of (\ref{A3v6}); then, the
optimal state variable $x(t)$ can be obtained from (\ref{xdotf}),
by
\begin{equation}
\dot{x}=f(x,u^{*}(x,\lambda^{*}(x,t),t),t).\label{xdotf2}
\end{equation}
Note that (\ref{xdotf2}) can be seen as the Pontryagin equation for $x(t)$, in the sense that $u^{*}$ and $\lambda^{*}$ are chosen in such a way that they are optimal closed-loop strategies, that is, these strategies maximize or minimize the action.\\
Thus, there exist three types of strategies in optimal control theory: open-loop ($x = x(t), \lambda = \lambda(t), u = u^*(t)$), inert closed-loop ($x = x(t), \lambda = \lambda(t), u = u^*(x(t), \lambda(t), t)$) and $\lambda$ closed-loop  ($x = x(t), \lambda = \lambda(x(t), t), u = u^*(x(t), \lambda(x(t), t), t)$) strategies. The first two are completely equivalent because these give the same dynamical equations.

\section{Dynamic Optimization: The Bellman approach}

A second approach to the optimization problem comes from dynamic programming theory, and was developed by Richard Bellman \cite{Bellman}. In this case, the fundamental variable is the optimal value of the action, defined by
\begin{equation}
J(x_{0},t_{0})=\max_{u}\Bigg(\int_{t_{0}}^{t}F(x,u,t)\ dt\Bigg),
\end{equation}
subject to (\ref{xdotf}), with initial condition $x(t_{0})=x_{0}$.\\
The optimality principle of Bellman implies that $J(x,t)$ satisfies the Hamilton--Jacobi--Bellman equation\cite{kamien}
\begin{equation}
\max_{u}\Bigl(F(x,u,t)+\frac{\partial J(x,t)}{\partial x}f(x,u,t)\Bigr)=-\frac{\partial J(x,t)}{\partial t}.\label{HJB}
\end{equation}
The left-hand side of Equation (\ref{HJB}) is just the maximization
of the Hamiltonian (\ref{pontryuv2}) with respect to the control
variable $u$, where the Lagrangian multiplier $\lambda$ of the Pontryagin approach must be identified with $\frac{\partial J(x,t)}{\partial x}$.
Thus, the Bellman theory is (from the Pontryagin perspective)
a model which has a closed-loop $\lambda$-strategy $\lambda(x,t)=\frac{\partial J(x,t)}{\partial x}$. \\
By maximizing and solving for the optimal control variable in the left-hand side of (\ref{HJB}) as $u^{*}=u^{*}(x,t)=u^{*}(x,\lambda(x,t),t)=u^{*}(x,\frac{\partial J(x,t)}{\partial x},t)$,
the Hamilton--Jacobi--Bellman equation is
\begin{equation}
F(x,u^{*},t)+\frac{\partial J(x,t)}{\partial x}f(x,u^{*},t)=-\frac{\partial J(x,t)}{\partial t}.\label{HJBv3}
\end{equation}
Differentiating in (\ref{HJBv3}) with respect to $x$, one gets \begin{small}
	\[
	\begin{array}{l}
	{\displaystyle \frac{\partial F(x,u^{*},t)}{\partial x}+{\displaystyle \frac{\partial^{2}J(x,t)}{\partial x^{2}}f(x,u^{*},t)+\frac{\partial J(x,t)}{\partial x}\frac{\partial f(x,u^{*},t)}{\partial x}}}\\
	\\
	+\Bigg({\displaystyle \frac{\partial F(x,u^{*},t)}{\partial u^{*}}+{\displaystyle \frac{\partial J(x,t)}{\partial x}\frac{\partial f(x,u^{*},t)}{\partial u^{*}}\Bigg)\frac{du^{*}(x,t)}{dx}}}\\
	\\
	=-{\displaystyle \frac{\partial^{2}J(x,t)}{\partial x\partial t}.}
	\end{array}
	\]
\end{small} Using the fact that $u^{*}$ is optimal and replacing
$\lambda(x,t)=\frac{\partial J(x,t)}{\partial x}$, one obtains
\[
{\displaystyle \frac{\partial F(x,u^{*},t)}{\partial x}+\frac{\partial\lambda(x,t)}{\partial x}f(x,u^{*},t)+\lambda\frac{\partial f(x,u^{*},t)}{\partial x}=-{\displaystyle \frac{\partial\lambda(x,t)}{\partial t},}}
\]
or
\[
\frac{dH^{*}(x,t)}{dx}+\frac{\partial\lambda(x,t)}{\partial t}=0.
\]
This equation is identical to (\ref{A3v6}). Thus, this optimal consistency
condition for the closed-loop $\lambda^{*}$-strategy is nothing more
than the derivative of the Hamilton--Jacobi--Bellman equation. Then, Equation (\ref{A3v6}) can be written, according to (\ref{HJBv3}),
as
\begin{equation}
\frac{d}{dx}\Bigg(F(x,u^{*},t)+\frac{\partial J(x,t)}{\partial x}f(x,u^{*},t)+\frac{\partial J(x,t)}{\partial t}\Bigg)=0.
\end{equation}
Integrating the above equation gives finally
\[
F(x,u^{*},t)+\frac{\partial J(x,t)}{\partial x}f(x,u^{*},t)+\frac{\partial J(x,t)}{\partial t}=g(t),
\]
where $g(t)$ is an arbitrary, time-dependent function. \\ 
Then, for an optimal closed-loop $\lambda^{*}$-strategy, the optimization problem gives a non-homogeneous Hamilton--Jacobi--Bellman equation (if one identifies the Lagrangian multiplier with $\frac{\partial J(x,t)}{\partial x}$). The Bellman maximum principle instead gives a homogeneous Hamilton--Jacobi--Bellman equation.

\section{Classical Mechanics and Open/Closed-Loop Strategies}

The classical equation of motion of a generic physical system can be obtained from the Hamiltonian action $A[x,p_{x}]$
defined by
\begin{equation}
A[x,p_{x}]=\int p_{x}\dot{x}-H(x,p_{x})\ dt,\label{Actionxp}
\end{equation}
(where $H$ is the Hamiltonian function) as an optimization problem. In this case, by expanding the Hamiltonian in a Taylor series, keeping the first-order terms only and integrating by parts, one has that the variation of the action is
\begin{equation}
\delta A={\displaystyle \int_{t_{0}}^{t_{1}}\Bigl[\left(\frac{\partial H}{\partial p_{x}}-\dot{x}\right)\delta p_{x}+\left(\frac{\partial H}{\partial x}+\dot{p_{x}}\right)\delta x\Bigr]dt+\bigl[p_{x}(t_{2})\delta x(t_{2})-p_{x}(t_{1})\delta x(t_{1})\bigr]}\label{deltaAdeltax}
\end{equation}
To maximize the action, all the first-order terms in $\delta x$ and $\delta p_{x}$ must vanish. For this, it is supposed that: 
\begin{itemize}
	\item the end points of the curve $x(t)$ are fixed, so $\delta x(t_{2})=0$
	and $\delta x(t_{1})=0$, and
	\item the variables $x$ and $p_{x}$ are considered independent,
\end{itemize}
so $\delta x$ and $\delta p_{x}$ are linearly independent, from
which the Hamiltonian equations of motion are obtained from (\ref{deltaAdeltax}) as:
\begin{equation} \label{hamiltoneq}
\dot{x}=\frac{\partial H(x,p_{x})}{\partial p_{x}}, \ \ \ \ \ \ \ \ \ \ \dot{p}_{x}=-\frac{\partial H(x,p_{x})}{\partial x}.
\end{equation}
From an economic or dynamic optimization point of view, the problem of optimizing the action (\ref{Actionxp}) is analogous to an optimal control problem, but without a control variable $u$. As we have seen in the previous section, the solutions of the control problem are given by the Pontryagin equations (\ref{pontry}) plus the transversality condition. Note that the first two Pontryagin equations in (\ref{pontry}) are precisely the Hamiltonian equations of motion (\ref{hamiltoneq}) if one identifies the Lagrangian multiplier $\lambda$ with the canonical momentum $p_{x}$. This implies that $\lambda$ corresponds to a $\lambda$
open-loop strategy for the system ($x(t)$ and $\lambda(t) = p_x(t)$ are independent),
which is consistent with the Hamiltonian equations of motions. Thus, the Hamiltonian theory is a open-loop model similar to Pontryagin's  theory.\\

So one may ask: is it possible that closed-loop strategies can occur in physical systems as they do in economic systems? \\
As is shown in this paper, the answer is affirmative, and they appear naturally in the context of canonical transformations. In order to give a first clue to the answer of the above question, suppose that one imposes a constraint over the phase space of the form
\begin{equation}
\Phi(x,p_{x},t)=0.\label{constraintpx}
\end{equation}
This constraint represents at each time a surface in the phase space
(a line in our bi-dimensional case) where the system can evolve (see
Figure 1).
\begin{figure}[h!]
	\centering \includegraphics[scale=0.9]{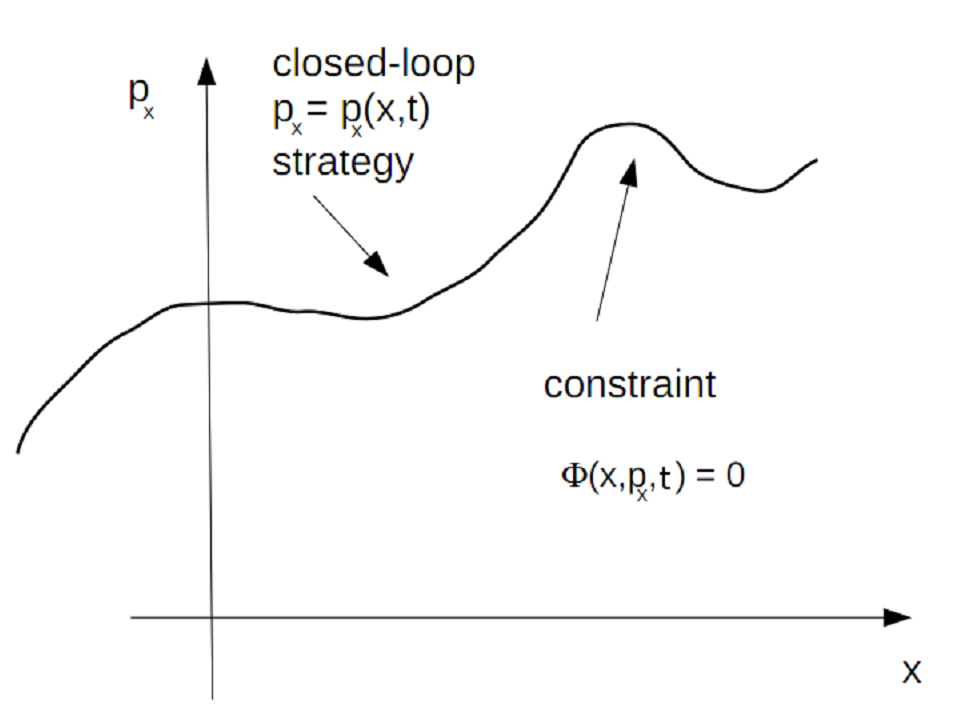}
	\caption{Line constraint over the phase space.}
	\label{fig:1}
\end{figure}
Actually one can write $p$ in terms of $x$ from (\ref{constraintpx})
by solving the constraint:
\begin{eqnarray}
p_{x}=p_{x}(x,t).\label{pxt}
\end{eqnarray}
Using the analogy of the momentum $p$ with the Lagrangian multiplier
$\lambda$, Equation (\ref{pxt}) corresponds to the closed-loop $\lambda$-strategy
from an economic point of view. Thus $x$ and $p$ are not independent
in this case and their variations are related by $\delta p_{x}=\frac{\partial p_{x}}{\partial x}\delta x$. By replacing in (\ref{deltaAdeltax}) and using the fact that the endpoints of
$x$ are fixed, one arrives at
\begin{equation}
\delta A={\displaystyle \int_{t_{0}}^{t_{1}}\Bigl[-\frac{\partial p_{x}}{\partial t}-\frac{\partial H}{\partial x}-\frac{\partial H}{\partial p_{x}}\frac{\partial p_{x}}{\partial x}\Bigr]\delta x\ dt.}\label{deltaAdeltaxclosed}
\end{equation}
By defining the reduced Hamiltonian $H^{*}(x,t)=H(x,p_{x}(x,t),t)$,
Equation (\ref{deltaAdeltaxclosed}) can be written as
\begin{equation}
\delta A={\displaystyle \int_{t_{0}}^{t_{1}}\Bigl[-\frac{\partial p_{x}}{\partial t}-\frac{\partial H^{*}}{\partial x}\Bigr]\delta x\ dt.}\label{deltaAdeltaxclosed 2}
\end{equation}
The optimization of the action then gives 
\begin{equation}
-\frac{\partial p_{x}}{\partial t}-\frac{\partial H^{*}}{\partial x}=0.\label{consistencypxt}
\end{equation}
This last equation is a consistency condition that the closed-loop strategy
(\ref{pxt}) must satisfy, to give an extremal of the action. That
is, if $p_{x}=p_{x}(x,t)$ satisfies (\ref{consistencypxt}), then
$p_{x}=p_{x}(x,t)$ is an optimal closed-loop strategy. \\
Now, if the closed-loop momentum strategy $p_{x}=p_{x}(x,t)$ is
just the derivative of some function $S(x,t)$, such as
\begin{eqnarray}
p_{x}(x,t)=\frac{\partial S(x,t)}{\partial x},\label{derivativecondition}
\end{eqnarray}
condition (\ref{consistencypxt}) gives
\begin{equation}
\frac{\partial}{\partial x}\Bigl[-\frac{\partial S(x,t)}{\partial t}-H^{*}(x,t)\Bigr]=0.
\end{equation}
So, one obtains, by integration, that
\begin{equation}
-\frac{\partial S(x,t)}{\partial t}-H(x,\frac{\partial S(x,t)}{\partial x})=g(t)\label{hamiltonjacobig(t)}
\end{equation}
for some function $g$ of time. The above equation is just an inhomogeneous Hamilton--Jacobi equation. Thus, the derivative of the Hamilton--Jacobi equation can be seen as the consistency condition to give to the action an extremal in the closed-loop $p_{x}$-strategy case. Also, this little analysis implies that closed-loop momentum strategies are closely related to the Hamilton--Jacobi equation.\\
Now, it is well known that the Hamilton--Jacobi equation appears in classical mechanics in the context of the canonical transformations, but how and where do the closed-loop $p_{x}=p_{x}(x,t)$ strategies appear there?
In the next section, a short review of canonical transformations will be given and it will be elucidated how the closed-loop strategies appear in that context.

\section{Canonical Transformations and Closed-Loop Strategies}

Consider a general coordinate transformation on the phase space
\begin{equation}
Q=Q(x,p_{x})\ \ \ \ \ P=P(x,p_{x}).\label{canonicaltransf}
\end{equation}
The transformation is called canonical if the Hamiltonian equations of motions are invariant under (\ref{canonicaltransf}), that is, if
\begin{equation}
\dot{Q}=\frac{\partial\tilde{H}}{\partial P}, \ \ \ \ \ \ \ \ \ \  \dot{P}=-\frac{\partial\tilde{H}}{\partial Q},
\end{equation}
where
\begin{equation}
\tilde{H}(P,Q,t)=H(x,p_{x},t)+\frac{\partial F(x,Q,t)}{\partial t}.
\end{equation}
The function $F$ is called the generator of the canonical transformation, and the coordinate transformation (\ref{canonicaltransf}) can be reconstructed from $F$ through Equations \cite{walecka} \  \cite{goldstein}
\begin{equation}
p_{x}=\frac{\partial F(x,Q,t)}{\partial x} , \ \ \ \ \ \ \ \ \ \ -P=\frac{\partial F(x,Q,t)}{\partial Q}.
\end{equation}
One must note at this point that the Hamiltonian equations (\ref{hamiltoneq}) refer to a unique coordinate system
(in this case, a Cartesian coordinate system). So, the Hamilton equations in (\ref{hamiltoneq}) are ``single observer'' equations. Instead, the canonical transformation brings a new second observer into the problem, because one has two different coordinate systems:
the initial Cartesian $(x,p_{x})$ and the second one $(Q,P)$. So
the theory of canonical transformations is a ``two observers'' view
of the classical mechanics. And this characteristic induces the closed-loop
$p_{x}$-strategies from a pure classical point of view (closed-loop
$p_{x}$-strategies can also be induced from Quantum Mechanics to
the classical realm, as we shall see later). \\
The Hamilton--Jacobi theory relies on the huge freedom that exists
in choosing $F(x,Q,t)$. In fact this theory doesn't work
directly with $F$ but with its Legendre transformation $S$ defined by
\begin{equation}
S(x,P,t)=F(x,Q,t)+PQ
\end{equation}
In this case, the canonical transformation is reconstructed via the equations
\begin{equation}
p_{x}=\frac{\partial S(x,P,t)}{\partial x}\label{reconstruc1}
\end{equation}
\begin{equation}
Q=\frac{\partial S(x,P,t)}{\partial P}\label{reconstruc2}
\end{equation}
and the respective Hamiltonians are related by
\begin{equation}
\tilde{H}(P,Q,t)=H(x,p_{x},t)+\frac{\partial S(x,P,t)}{\partial t}\label{Hamiltoniantilde}
\end{equation}
Note again that one has huge freedom in choosing $S$. The Hamilton--Jacobi theory corresponds to the choice of $S$
that makes the second-observer Hamiltonian $\tilde{H}(Q,P)$
equal zero:
\begin{equation}
\tilde{H}(P,Q,t)=0.\label{Hamiltoniantildeequal0}
\end{equation}
So the equations of motion for the second observer are
\begin{equation}
\dot{Q}=0 ,  \ \ \ \ \ \ \ \ \ \ \  \ \dot{P}=0.
\end{equation}
Thus, for the second observer, the dynamical variable remain constant in time:
\begin{equation} \label{Qequal0}
Q(t)=Q_{0},  \ \ \ \ \ \ \ \ \ \ \  \ P(t)=P_{0}.
\end{equation}
But what does the first observer see? First, due to Equations (\ref{Qequal0}), the coordinate transformations (\ref{canonicaltransf})
give
\begin{equation} \label{closeloop1}
Q_{0}=Q(x,p,t),  \ \ \ \ \ \ \ \ \ \ \  \ P_{0}=P(x,p,t),
\end{equation}
but each of these equations defines constant coordinate lines. These are constraints over the phase space of the first observer, from which one can generate two different closed-loop $p_{x}$-strategies
according to Equation (\ref{constraintpx})! (See Figure 2).
\begin{figure}[h]
	\centering \includegraphics[scale=0.9]{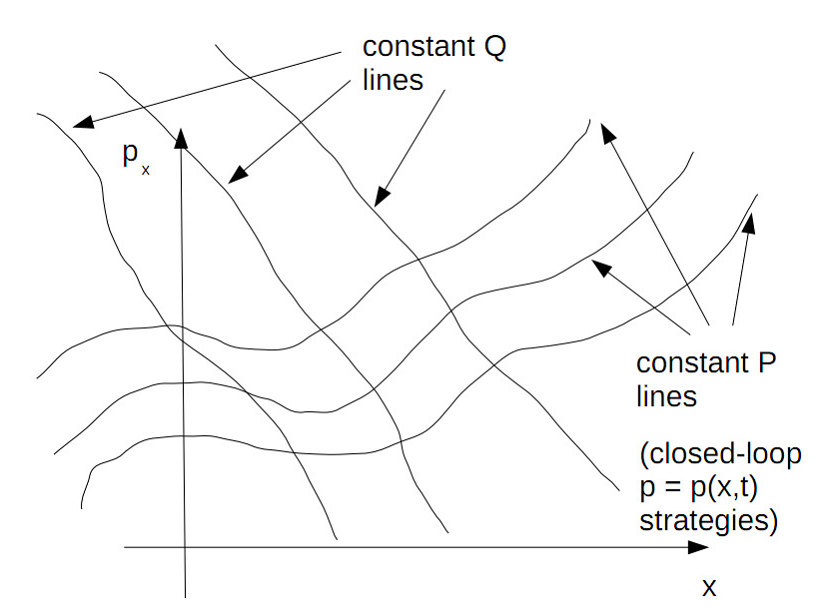}
	\caption{Constraints over the phase space.}
	\label{fig:2}
\end{figure}
Thus, the ``two observers'' perspective of classical mechanics,
through the method of canonical transformation, is responsible for
the generation of the closed-loop $p_{x}$-strategies. \\
From (\ref{closeloop1}) it is not clear if
the $p_{x}$ closed-loops strategies thus generated satisfy the consistency
condition (\ref{consistencypxt}) or if they satisfy the second condition
(\ref{derivativecondition}) to obtain a Hamilton--Jacobi equation
as in (\ref{hamiltonjacobig(t)}) for $S$. Instead, one can see these constant coordinate lines in terms of Equations
(\ref{reconstruc1}) and (\ref{reconstruc2}). In fact, these equations
are equivalent to (\ref{closeloop1}), because
the canonical transformation can be reconstructed from (\ref{reconstruc1})
and (\ref{reconstruc2}). \\
A constant $P_{0}$-line in (\ref{closeloop1}) is equivalent (from
(\ref{reconstruc1})) to
\begin{equation} \label{pxconsistentstrategy2}
p_{x}=\frac{\partial S(x,P_{0},t)}{\partial x},
\end{equation}
thus, the constant $P_{0}$-line in (\ref{closeloop1}) satisfies
(\ref{derivativecondition}). From (\ref{Hamiltoniantilde}) and (\ref{Hamiltoniantildeequal0}),
$S(x,t)$ satisfies the Hamilton--Jacobi equation as in
(\ref{hamiltonjacobig(t)}) with $g(t)=0$. This implies that a closed-loop strategy generated by (\ref{pxconsistentstrategy2}) satisfies the consistency relation (\ref{consistencypxt}), so (\ref{pxconsistentstrategy2}) defines a true optimal closed-loop $p_{x}$-strategy.\\
For the $Q_{0} $ constant coordinate line in (\ref{closeloop1}) one has
however due to (\ref{reconstruc2}) that
\begin{equation}
Q_{0}=\frac{\partial S(x,P,t)}{\partial P},\label{reconstruc22}
\end{equation}
but from this equation one cannot obtain $p_{x}$ in terms of $x$, so this line does not generates a $p_{x}$ closed-loop strategy at all. This is due to the structure of Equations (\ref{reconstruc1}) and (\ref{reconstruc2}). In order to reconstruct the canonical transformation
(\ref{canonicaltransf}) it is required to invert the system (\ref{reconstruc1}),
(\ref{reconstruc2}). Note that only from (\ref{reconstruc1}) can the
momentum $P$ be written in terms of the first-observer variables
as $P=P(x,p_{x},t)$. But from (\ref{reconstruc2}) alone one cannot
solve $Q$ in terms of the $x,p_{x}$. The other equation (\ref{reconstruc1})
is needed to do that. Thus, a constant $Q_{0}$-line alone can not generate a true closed-loop $p_{x}=p_{x}(x,t)$ strategy. \\

In this way, closed-loop $p_x$-strategies appear in classical mechanics as a consequence of the two observers interpretation of the canonical transformation theory. These closed-loop strategies are inert in the same way that the optimal $u^*$ closed-loop ones are inert in control theory. This is because both closed-loop approaches $p_x$ and $u^*$ gives the same dynamical equation of the open-loop case. For the closed-loop $p_x$ case, the open-loop dynamics analogous to the given by Pontryagin's equations are just provided by the Hamiltonian equations of motion of the first observer in the $(x, p_x)$ phase space.

\section{Quantum Mechanics and Closed-Loop Strategies}

In this section, the origins of the Hamilton--Jacobi--Bellman equation that appears in the limit $\hbar\rightarrow0$ in the quantum phenomena,
will be explained as a consequence of the emergence of closed-loop
$p_{x}$-strategies in the quantum world. \\
Consider the Schrödinger equation for a non-relativistic particle
of mass $m$:
\begin{equation}
-\frac{\hbar^{2}}{2m}\frac{\partial^{2}\Psi(x,t)}{\partial x^{2}}+U(x)\Psi(x,t)=i\hbar\frac{\partial\Psi(x,t)}{\partial t} .\label{Schrodinger}
\end{equation}
Writing the wave function in the form
\begin{equation}
\Psi(x,t)=e^{\frac{i}{\hbar}S(x,t)},\label{Sdext}
\end{equation}
and by substituting (\ref{Sdext}) into the Schrödinger equation, the
following equation Quantum Hamilton--Jacobi equation for $S(x,t)$
is obtained:
\begin{equation}
\frac{1}{2m}\Bigl(\frac{\partial S(x,t)}{\partial x}\Bigr)^{2}+U(x)-\frac{i\hbar}{2m}\frac{\partial^{2}S(x,t)}{\partial x^{2}}=-\frac{\partial S(x,t)}{\partial t} .\label{quantumHJB}
\end{equation}
Note that this equation is completely equivalent to Schrödinger's equation,
but here, the classical and quantum realms can be clearly identified.
In fact, by taking the limit $\hbar\rightarrow0$ in (\ref{quantumHJB}),
one gets
\begin{equation}
\frac{1}{2m}\Bigl(\frac{\partial S(x,t)}{\partial x}\Bigr)^{2}+U(x)=-\frac{\partial S(x,t)}{\partial t} . \label{classicalHJB}
\end{equation}
Equation (\ref{classicalHJB}) is just the classical Hamilton--Jacobi equation
\begin{equation}
H(x,\frac{\partial S(x,t)}{\partial x})=-\frac{\partial S(x,t)}{\partial t} , \label{classicalHJBgeneral}
\end{equation}
for the classical Hamiltonian function associated to the non-relativistic
particle
\begin{equation}
H(x,p_{x})=\frac{p_{x}^{2}}{2m}+U(x)  , \label{hamiltonianfunction}
\end{equation}
where one must identify $p_{x}$ with the derivative of $S$
\begin{equation}
p_{x}=p_{x}(x,t)=\frac{\partial S(x,t)}{\partial x}.\label{identification}
\end{equation}
to make contact with the classical  Hamiltonian theory. And it is precisely this identification which generates the closed-loop
$p_{x}$-strategy through (\ref{identification}). Note that it is induced from the quantum realm to the classical world, in the
limit $\hbar\rightarrow0$, through the classical Hamilton--Jacobi
equation. The identification in (\ref{identification}) is a pure quantum phenomenon. In fact, considering the momentum operator 
\begin{equation}
\hat{p}_{x}=-i\hbar\frac{\partial}{\partial x} ,
\end{equation}
this operator is characterized by its eigenfunctions and eigenvalues:
\begin{equation}
\hat{p}_{x}\Phi_{p_{x}}(x)=p_{x}\Phi_{p_{x}}(x).
\end{equation}
where the solution of this equation gives
\begin{equation}
\Phi_{p_{x}}(x)=e^{\frac{i}{\hbar}p_{x}x}.\label{momentumeigenfunction}
\end{equation}
In this context, $p_{x}$ and $x$ are independent variables. In fact
the eigenfunction (\ref{momentumeigenfunction}) corresponds to states with well
defined values of the momentum. \\
Note now that if one applies the momentum operator to a generic wave function $\Psi$ which is a solution of the Schrödinger equation (written in the ``momentum
form'' (\ref{Sdext})), one obtains
\begin{equation}
\hat{p}_{x}\Psi(x,t)=\frac{\partial S(x,t)}{\partial x}\Psi(x,t).\label{diagonalpPsi}
\end{equation}
By looking at the wave function as a vector with a continuous index $x$, the above equation implies that (locally at each point $x$) the momentum operator is diagonal, so that any wave function can be seen as an eigenstate of the momentum operator with momentum eigenvalue $\frac{\partial S(x,t)}{\partial x}$. Thus, one must identify the momentum eigenvalue $p_{x}$ in this quantum
state with the derivative of the $S$ function through (\ref{identification}).
It is just this identification which generates the closed-loop $p_{x}$-strategies directly
in the quantum world. \\
On the other hand, the same Heisenberg canonical commutation relations
\begin{equation}
[\check{x},\check{P_{x}}]=\check{x}\check{P_{x}}-\check{P_{x}}\check{x}=i\hbar\ \check{I}.\label{HCR}
\end{equation}
can be seen as a constraint in the non-commutative phase space $(\check{x},\check{P_{x}})$.
Thus, from (\ref{HCR}) one could ``solve'' the momentum operator
$\check{P_{x}}$ in terms of the $\check{x}$ operator. This necessarily
implies the existence of a certain relation between $\check{P_{x}}$
and $\check{x}$ or between their eigenvalues. The representation of
the canonical operator as a differential operator acting on a function
space or Hilbert space as
\begin{equation}
\check{x}\rightarrow x,\ \ \ \ \ \hat{p}_{x}\rightarrow-i\hbar\frac{\partial}{\partial x}
\end{equation}
is equivalent to solving the constraint (\ref{HCR}), because on
any wave function $\Psi(x,t)$, Equation (\ref{HCR}) is satisfied
identically. The memory of the quantum constraint (\ref{HCR}) is then transferred in a local way to the momentum eigenvalue, according to (\ref{diagonalpPsi}).
In a sense, the representation of the wave function as $\Psi(x,t)=e^{\frac{i}{\hbar}S(x,t)}$
locally diagonalizes the momentum operator $\check{p}_{x}$ over any quantum
state, and (\ref{quantumHJB}) is just the Schrödinger equation in
this diagonal basis. Note that all this is a kinematic effect created by the Heisenberg commutation relation (\ref{HCR}); the
dynamical effects appear when the explicit form of $S(x,t)$ is
needed, and for that, one must solve the full Quantum Hamilton--Jacobi equation
(\ref{quantumHJB}) explicitly.  Note that Quantum Mechanics, as in (\ref{HCR}), can be viewed as a constrained system in a non-commutative space, so, one would apply a generalization of Dirac's method \cite{dirac1} \ \cite{dirac2} \ \cite{teitelboim} \ \cite{rothe} to non-commutative spaces \cite{contreras3} to study quantum mechanical systems.\\

We can say, then, that closed-loop $p_{x}$-strategies correspond to a 
pure quantum phenomenon and are a consequence of Heisenberg's uncertainty
principle. In an arbitrary quantum state, momentum and position cannot
be independent: they are related through the non-commutative character
of the position and momentum operators. In a more defined momentum
state, a less defined position state would emerge. Thus, these two
variables must depend on one another in some way. Relation (\ref{identification})
is tantamount to  a conversation between them. Only in a pure-momentum
state, as given in (\ref{momentumeigenfunction}), does the link disappear
and position and momentum become independent variables. In fact, in
a pure-momentum state, $p_{x}(x,t)=p_{x}^{0}$ is constant, that is:
all the eigenvalues are the same, the value $x$ of the position (so
its matrix is a multiple of the identity matrix), so Equation (\ref{identification})
gives
\begin{eqnarray}
S(x,t)=p_{x}^{0}x+\phi(t)
\end{eqnarray}
as a solution, where $\phi(t)$ is some function of time. Thus, the
wave function is
\begin{equation}
\Psi(x,t)=e^{\frac{i}{\hbar}(p_{x}^{0}x+\phi(t))}=e^{\frac{i}{\hbar}\phi(t)}\Phi_{p_{x}^{0}}(x)
\end{equation}
which is the same momentum eigenstate amplified by a temporal arbitrary
phase. Then the linear character of $S(x,t)$ in terms of $x$ implies
that $p_{x}$ and $x$ are independent variables, and no closed-loop
$p_{x}$-strategy exists in this case. \\
The same can be said for a pure-position eigenstate. Thus, closed-loop
$p_{x}$-strategies are an inherent part of the quantum mechanical
world and permeate the classical world in the limit $\hbar\rightarrow0$
through the Hamilton--Jacobi equation. In the following two sections, we analyze some common textbook examples from closed-loop strategies' point of view.

\section{The Stationary Case}

The quantum and classical Hamilton--Jacobi equations (\ref{quantumHJB})
and (\ref{classicalHJB}) are non-stationary equations, that is, they
depend explicitly on time. In Quantum Mechanics, stationary
states play a fundamental role. They are defined by
\begin{equation}
\Psi(x,t)=e^{-\frac{i}{\hbar}Et}\Phi(x).\label{eEt}
\end{equation}
By substituting this into the time-dependent Schrödinger equation (\ref{Schrodinger}),
the time independent or stationary Schrödinger equation is obtained:
\begin{equation}
-\frac{\hbar^{2}}{2m}\frac{\partial^{2}\Psi(x)}{\partial x^{2}}+U(x)\Psi(x)=E\Psi(x).\label{Schrodingernotime}
\end{equation}
Now by writing
\begin{equation}
\Psi(x)=e^{\frac{i}{\hbar}W(x)}\label{ewdex}
\end{equation}
and substituting into (\ref{Schrodingernotime}) the stationary Quantum
Hamilton--Jacobi equation holds:
\begin{equation}
\frac{1}{2m}\Bigl(\frac{\partial W(x)}{\partial x}\Bigr)^{2}+U(x)-\frac{i\hbar}{2m}\frac{\partial^{2}W(x)}{\partial x^{2}}=E.\label{quantumHJBnotime}
\end{equation}
Taking again the classical limit $\hbar\rightarrow0$ in (\ref{quantumHJBnotime})
the stationary Classical Hamilton--Jacobi equation appears:
\begin{equation}
\frac{1}{2m}\Bigl(\frac{\partial W(x)}{\partial x}\Bigr)^{2}+U(x)=E.\label{classicalHJBnotime}
\end{equation}
Due to (\ref{eEt}) and (\ref{ewdex}), we have $\Psi(x,t)=e^{\frac{i}{\hbar}S(x,t)}=e^{\frac{i}{\hbar}(W(x)-Et)}$
which implies that
\begin{equation}
S(x,t)=W(x)-Et
\end{equation}
for the stationary case. In this case, the closed-loop $p_{x}$-strategies
are given by $p_{x}= p_x(x,t) =\frac{\partial W(x)}{\partial x}$.

\section{The Non-Stationary Case}

But what about the non-stationary closed-loop $p_{x}$-strategies
in the classical limit? In order to analyze this case, consider the
example of a free particle, that is, $U(x)=0$. The non-stationary
classical Hamilton--Jacobi equation is now
\begin{equation}
\frac{1}{2m}\Bigl(\frac{\partial S(x,t)}{\partial x}\Bigr)^{2}=-\frac{\partial S(x,t)}{\partial t}.\label{freeCHJ}
\end{equation}
One can find a solution of the form \ $S(x,t)=\frac{1}{2}a(t)x^{2}$, \ so by substituting in (\ref{freeCHJ}) one gets $a(t)=\frac{-1}{(P_{0}-\frac{t}{m})}$,
so $S(x,t)=\frac{1}{2}\frac{-x^{2}}{(P_{0}-\frac{t}{m})}$, and the corresponding closed-loop $p_{x}$-strategy is
\begin{equation}
p_{x}(x,t) = \frac{\partial S(x,t)}{\partial x}=\frac{-x}{(P_{0}-\frac{t}{m})} .\label{pxdex}
\end{equation}
One can evaluate $x(t)$ using Equation (\ref{reconstruc22})
\begin{equation}
Q_{0}=\frac{\partial S(x,P_{0},t)}{\partial P_{0}} ,
\end{equation}
where the integration constant $P_{0}$ must be identified with the constant momentum for the second observer in the coordinate system $(Q,P)$. Thus
\begin{equation}
Q_{0}=\frac{1}{2}\frac{x^{2}}{(P_{0}-\frac{t}{m})^{2}} ,
\end{equation}
from which $x(t)$ is computed as
\begin{equation} \label{xdet}
x(t)=\sqrt{2\left(P_{0}-\frac{t}{m}\right)^{2}Q_{0}}=\sqrt{2Q_{0}}\left(P_{0}-\frac{t}{m}\right).
\end{equation}
The associated open-loop $p_{x}(t)$ strategy is found by $p_{x}(t) = p_x(x(t),t)$ similarly to Equation (\ref{openloop_u*(t)}). Thus, by substituting (\ref{xdet}) into (\ref{pxdex}):
\begin{equation}
p_{x}(t)=\frac{-x(t)}{(P_{0}-\frac{t}{m})}=\frac{-\sqrt{2Q_{0}}(P_{0}-\frac{t}{m})}{(P_{0}-\frac{t}{m})}=-\sqrt{2Q_{0}}=p_{0} ,
\end{equation}
so
\begin{equation} 
x(t)=-p_{0}\left(P_{0}-\frac{t}{m}\right)=p_{0}\frac{t}{m}+x_{0}, 
\end{equation}
where $x_{0}=-P_{0}p_{0}$. These last equations are the solutions
for the motion of a free particle of course! \\
Now from the Hamiltonian equations, one gets the open-loop dynamics for the free particle:
\begin{equation}
\dot{x}=\frac{p_{x}}{m},  \ \ \ \ \ \ \ \ \ \dot{p}_{x}=-\frac{\partial U(x)}{\partial x}=0 ,
\end{equation}
so the open-loops dynamics is
\begin{equation}
x(t)=\frac{p_{0}}{m}t+x_{0} ,  \ \ \ \ \ \ \ \ \  p_{x}(t) = p_{0} .
\end{equation}
Then, the open-loop $p_{x}$-strategy for the free particle coming
from the Hamiltonian equations of motions is identical to the non-stationary
closed-loop $p_{x}$-strategy coming from the non-stationary classical
Hamilton--Jacobi equation. In fact, any solution of the Classical Hamilton--Jacobi
equation would be an origin for a (stationary and non-stationary) closed-loop
$p_{x}$-strategy, and this one has to be equivalent with the open-loop
strategy coming from the Hamiltonian equations of motion. This is
because the Classical Hamilton--Jacobi equation corresponds to a ``two
observer'' point of view of classical mechanics. The function $S(x,t)$
is just the generator of the canonical transformation which
leaves the Hamiltonian equations invariant. Thus, the two schemes,
\begin{enumerate}
	\item the Hamiltonian ``one observer'' approach with its open-loop
	$p_{x}$-strategies and
	\item the Classical Hamilton--Jacobi ``two observer'' approach with its closed-loop $p_{x}$-strategies
\end{enumerate}
are equivalent, because they have the same equation of motion constructed
through the canonical transformation. Closed-loop strategies coming from the Hamilton--Jacobi equation are similar in character to the inert optimal $u^{*}(x,t)$ closed-loop strategies of the Pontryagin approach to optimal control theory, because they are equivalent to the open-loop ones.

\section{The Pure Quantum Limit and Closed-Loop Strategies}

In the previous section the classical limit $\hbar\rightarrow0$ was
taken and its characteristics were explored in terms of the closed-loop
strategies. In this section, the inverse limit is going to be taken,
that is: $\hbar>>1$, and the consequences of a higher non-commutative
system of quantum variables explored:
\begin{equation}
[\check{x},\check{P_{x}}]=i\hbar\ \check{I}.
\end{equation}

Consider the non-stationary Quantum Hamilton--Jacobi equation (\ref{quantumHJB}). Taking the limit $\hbar>>1$ and supposing that the time derivative
of the $S$ function has a higher value,
\begin{equation}
-\frac{i\hbar}{2m}\frac{\partial^{2}S(x,t)}{\partial x^{2}}=-\frac{\partial S(x,t)}{\partial t}.\label{SQHJ}
\end{equation}
or, what is the same,
\begin{equation}
-\frac{\hbar^{2}}{2m}\frac{\partial^{2}S(x,t)}{\partial x^{2}}=i\hbar\frac{\partial S(x,t)}{\partial t}.
\end{equation}
But this is the free-particle Schrödinger equation for $S(x,t)$!
In this way one can write
\begin{equation}
S(x,t)=e^{\frac{i}{\hbar}T(x,t)}
\end{equation}
where $T(x,t)$ satisfies the Quantum Hamilton--Jacobi equation
\begin{equation}
\frac{1}{2m}\Bigl(\frac{\partial T(x,t)}{\partial x}\Bigr)^{2}-\frac{i\hbar}{2m}\frac{\partial^{2}T(x,t)}{\partial x^{2}}=-\frac{\partial T(x,t)}{\partial t}.\label{TQHJ}
\end{equation}
and the wave function is given by
\begin{equation}
\Psi(x,t)=e^{\frac{i}{\hbar}S(x,t)}=e^{\frac{i}{\hbar}e^{\frac{i}{\hbar}T(x,t)}}
\end{equation}
The corresponding closed-loop $p_{x}$-strategy is
\begin{equation}
p_{x}=\frac{\partial S(x,t)}{\partial x}=\frac{i}{\hbar}\frac{\partial T(x,t)}{\partial x}S(x,t)
\end{equation}
Note that $\frac{\partial T(x,t)}{\partial x}$ can be interpreted
as a closed-loop $p_{T}$-strategy for the Quantum Hamilton--Jacobi
equation (\ref{TQHJ}). Thus, denoting the closed-loop $p_{x}$-strategy
for (\ref{SQHJ}) by $p_{S}(x,t)$, then both strategies are related
by
\begin{equation}
p_{S}(x,t)=\frac{i}{\hbar}\ S(x,t)\ p_{T}(x,t)\label{pSpTestrategy}
\end{equation}

Again, if the time derivative of $T$ has a higher value,
and as $\hbar>>1$, the Quantum Hamilton--Jacobi for $T(x,t)$ (\ref{TQHJ})
is in this limit again a Schrödinger equation,
\begin{equation}
-\frac{\hbar^{2}}{2m}\frac{\partial^{2}T(x,t)}{\partial x^{2}}=i\hbar\frac{\partial T(x,t)}{\partial t}.
\end{equation}
Hence, we can write
\begin{equation}
T(x,t)=e^{\frac{i}{\hbar}U(x,t)}
\end{equation}
where $U(x,t)$ satisfies
\begin{equation}
\frac{1}{2m}\Bigl(\frac{\partial U(x,t)}{\partial x}\Bigr)^{2}-\frac{i\hbar}{2m}\frac{\partial^{2}U(x,t)}{\partial x^{2}}=-\frac{\partial U(x,t)}{\partial t}.\label{UQHJ}
\end{equation}
and the wave function is
\begin{equation}
\Psi(x,t)=e^{\frac{i}{\hbar}S(x,t)}=e^{\frac{i}{\hbar}e^{\frac{i}{\hbar}e^{\frac{i}{\hbar}U(x,t)}}}\label{pSpU1}
\end{equation}
Putting $p_{U}=\frac{\partial U(x,t)}{\partial x}$ the closed-loop
strategy associated to the Quantum Hamilton--Jacobi equation for $U$
(\ref{UQHJ}), the corresponding closed-loop $p_{S}$-strategy is
then in this case
\begin{equation}
p_{S}(x,t)=\frac{i}{\hbar}p_{U}(x,t)S(x,t)T(x,t).\label{pSpU}
\end{equation}
Since one can keep iterating this procedure to infinity, quantum
mechanical systems can admit multistage closed-loop strategies and
they are connected in strongly non-linear way, as in (\ref{pSpU1}).

\section{Conclusions}

In this article, we developed an optimal control perspective on the dynamical behavior of classical and quantum physical systems. The most crucial element of this view is the presence of feed-backs characterized by open or closed-loop strategies in the system. \\ \\
Thus, in quantum theory, the closed-loop strategies appear naturally due thinking of Heisenberg's commutation relations as a constraint in a non-commutative phase space, so this implies that there is a relation between the momentum and the particle position for any quantum state. \\ \\
By taking the classical limit $\hbar\rightarrow0$ in the full Quantum Hamilton--Jacobi equation, one arrives at a closed-loop dynamics associated with the Classical Hamilton--Jacobi theory. The non-commutative character of quantum theory is transferred to the classical theory through the closed-loop $p_{x} = \frac{\partial S(x,t)}{\partial x}$-strategy. Since $S(x,t)$ satisfies the Classical Hamilton--Jacobi equation, the dynamics generated by $S(x,t)$ (by virtue of the properties of canonical transformations, whose generator is just $S(x,t)$) is completely equivalent to those open-loop dynamics dictated by the Hamiltonian equations of motion. \\ \\
From a purely classical point of view, the presence of these closed-loop strategies can be explained by the ``two observers'' character of the Hamilton--Jacobi theory. If the solutions of the equations of motion are constant for the first observer, then for the second one, their solutions look like constraints. This necessarily relates the momentum of the particle with its position for the second observer, generating in this way the closed-loop $p_{x}=\frac{\partial S(x,t)}{\partial x}$-strategy.


\begin{thebibliography}{99}
	
	\bibitem{stanley} Rosario N. Mantegna and H. Eugene Stanley, \textit{An Introduction to Econophysics} (Cambridge University Press, 2007).
	
	\bibitem{boucheaud} Jean-Philippe Boucheaud and Marc Potters, \textit{Theory
		of Financial Risk and Derivative Pricing: From Statistical Physics
		to Risk Management} (Cambridge University Press, 2009).
	
	\bibitem{baaquie} B. E. Baaquie, \textit{Quantum Finance: Path Integrals
		and Hamiltonians for Option and Interest Rates} (Cambridge University
	Press, 2007).
	
	\bibitem{kamien} Morton I. Kamien and Nancy L. Schwartz, \textit{The
		Calculus of Variations and Optimal Control in Economics and Management} (Dover, 2012).
	
	\bibitem{sheti} Suresh P. Sethi, \textit{Optimal Control Theory:
		Applications to Management Science and Economics} (Springer, 2009).
	
	\bibitem{caputo} Michael R. Caputo, \textit{Foundations of Dynamic
		Economic Analysis: Optimal Control Theory and Applications} (Cambridge University Press, 2005).
	
	\bibitem{weitzman} Martin L. Weitzman , \textit{Income, Wealth, and
		the Maximum Principle} (Harvard University Press, 2007).
	
	\bibitem{dockner} Engelbert J. Dockner, Steffen Jorgensen, and Ngo Van Long, \textit{Differential Games in Economics and Management
		Science} (Cambridge University Press, 2001).
	
	\bibitem{contrerasRPMV} M. Contreras, R. Pellicer, and M. Villena,
	``Dynamic optimization and its relation to classical and quantum
	constrained systems,'' Physica A \textbf{479}, 12--25 (2017).
	
	\bibitem{hojman}  Sergio A. Hojman, 
	``Optimal Control and Dirac's Theory of Singular Hamiltonian Systems,'' unpublished.
	
	\bibitem{teturo} Teturo Itami, 
	``Quantum Mechanical Theory of Nonlinear Control'' in \textit{IFAC Nonlinear Control Systems} (IFAC Publications, St. Petersburg, Russia, 2001), p. 1411.
	
	\bibitem{contreras2} M. Contreras and J. P. Peña,
	``The quantum dark side of the optimal control theory,'' Physica A \textbf{515}, 450--473 (2019).
	
	\bibitem{Pontryagin} L. S. Pontryagin, V. G. Boltyanskii, R. V. Gamkrelidze, and E. F. Mishchenko, \textit{The Mathematical Theory of Optimal Processes}, (CRC Press, 1987).
	
	\bibitem{Bellman} R. Bellman, ``The theory of dynamic programming,'' Bull. Am. Math. Soc. \textbf{60} (6), 503--516 (1954).
	
	\bibitem{erikson1} G. M. Erickson, ``Differential game models
	of advertising competitions,'' J. of Political Economy \textbf{8} (31), 637--654 (1973).
	
	\bibitem{walecka} Alexander L. Fetter and John Dirk Walecka, \textit{Theoretical
		Mechanics of Particles and Continua} (Dover, 2003).
	
	\bibitem{goldstein} Herbert Goldstein, \textit{Classical Mechanics}, 3rd ed. (Pearson, 2001).
	
	\bibitem{dirac1} P. A. M. Dirac,
	``Generalized Hamiltonian dynamics,'' Proc. Roy. Soc. London, A \textbf{246}, 326 (1958).
	
	\bibitem{dirac2} P. A. M. Dirac, \textit{Lectures on Quantum Mechanics} (Yeshiva University Press, New York, 1967).
	
	\bibitem{teitelboim} C. Teitelboim and M. Henneaux, \textit{Quantization of Gauge Systems} (Princeton University Press, 1994).
	
	\bibitem{rothe} H. J. Rothe  and K. D. Rothe, \textit{Classical and Quantum Dynamics of Constrained
		Hamiltonian Systems, World Scientific Lectures Notes in Physics}, vol. 81, (World Scientific, 2010).
	
	\bibitem{contreras3} M. Contreras G.,
	``Dirac's method in a non-commutative phase space,'' in preparation.		
	
\end{thebibliography}
\end{document}